\documentclass[12pt]{article}
\usepackage{xspace,vmargin,float,hhline}
\usepackage{amsmath,amssymb,latexsym,righttag,theorem}

\setpapersize{USletter}
\setmarginsrb{0.75in}{0.9in}{0.75in}{0.9in}{0pt}{0mm}{0pt}{10mm}

\newcommand{\rreal}     {{\mathbb{R}}}
\newcommand{\nats}      {{\mathbb{N}}}

\newcommand{\eg}{e.g.,\xspace}
\newcommand{\ie}{i.e.,\xspace}

\newcommand{\be}{\begin{equation}}
\newcommand{\ee}{\end{equation}}
\newcommand{\p}{\partial}

\newcommand{\Lnorm}{ {\mathcal L} }         

\theoremstyle{break}            
\theoremstyle{plain}
{\theorembodyfont{\sffamily\small}      }
{\theorembodyfont{\rmfamily}    }
\theorembodyfont{\itshape}      
\theoremstyle{marginbreak}      
\theoremstyle{change}
\theoremstyle{break}            

\theoremheaderfont{\scshape}

\newcommand{\nohyphens}{\hyphenpenalty=10000\exhyphenpenalty=10000\relax}

\renewcommand{\appendix}{%
   \renewcommand{\section}{%
        \newpage\thispagestyle{plain}%
        \secdef\Appendix\sAppendix}%
   \setcounter{section}{0}%
   \renewcommand{\thesection}{\Alph{section}}%
}

\newcommand{\Appendix}[2][?]{%
   \refstepcounter{section}%
   \addcontentsline{toc}{appendix}%
        {\protect\numberline{\appendixname~\thesection} #1}%
   {\flushright\large\bfseries\appendixname\ \thesection\par
    \nohyphens\centering#2\par}%
   \sectionmark{#1}\vspace{\baselineskip}
}

\newcommand{\sAppendix}[1]{%
        {\flushright\large\bfseries\appendixname\par
         \nohyphens\centering#1\par}%
        \vspace{\baselineskip}
}

\newcounter{theoremm}
\newcounter{remarkk}
\newcounter{lemmaa}

\newcounter{examplee}
\newcounter{applicationn}
\newcounter{definitionn}
\newcounter{exercisee}
\newcounter{assumee}

\def\thetheoremm{\arabic{theoremm}}
\def\theremarkk{\arabic{remarkk}}
\def\thelemmaa{\arabic{lemmaa}}
\def\thecorollaryy{\arabic{corollaryy}}
\def\theexamplee{\arabic{examplee}}
\def\theapplicationn{\arabic{applicationn}}
\def\thedefinitionn{\arabic{definitionn}}

\def\theassumee{\arabic{assumee}}

\newenvironment{theoremm}{\refstepcounter{theoremm}   \vspace{2ex} {\bf Theorem \thetheoremm}: \em }{\vspace{1ex} }

\newenvironment{remarkk}{\refstepcounter{remarkk}   \vspace{2ex} {\bf Remark \theremarkk}: }{\vspace{1ex} }

\newenvironment{proof}{ \par\noindent \space \space {\it Proof: }}
{}

\begin{document}

\title{Adaptive Control for a Class of Nonlinear Systems with a Time-Varying
  Structure}

\makeatletter

\author{Ra\'ul Ord\'o\~nez\thanks{R. Ord\'o\~nez is with the
    Dept. of Electrical and Computer Engineering, Rowan
    University, 201 Mullica Hill Road, Glassboro, New Jersey, 08028
    (ordonez@rowan.edu).} \and Kevin M. Passino\thanks{K. M. Passino
    is with the Dept. of Electrical Engineering, The Ohio State
    University, 2015 Neil Ave., Columbus, OH 43210
    (k.passino@osu.edu).}}

\makeatother

\date{}

\maketitle

\begin{abstract}
  In this paper we present a direct adaptive control method for a class of
  uncertain nonlinear systems with a time-varying structure.  We view the
  nonlinear systems as composed of a finite number of ``pieces,'' which are
  interpolated by functions that depend on a possibly exogenous scheduling
  variable. We assume that each piece is in strict feedback form, and show
  that the method yields stability of all signals in the closed-loop, as well
  as convergence of the state vector to a residual set around the equilibrium,
  whose size can be set by the choice of several design parameters. The class
  of systems considered here is a generalization of the class of strict
  feedback systems traditionally considered in the backstepping literature. We
  also provide design guidelines based on $\Lnorm_\infty$ bounds on the
  transient.
\end{abstract}

\section{Introduction}

The field of nonlinear adaptive control developed rapidly in the last
decade. The paper \cite{Polycarpou91} and others gave birth to an
important branch of adaptive control theory, the nonlinear on-line
function approximation based control, which includes neural (\eg in
\cite{P96}) and fuzzy (\eg in \cite{Su94}) approaches (note that there
are several other relevant works on neural and fuzzy control, many of
them cited in the references within the above papers). The neural and
fuzzy approaches are most of the time equivalent, differing between
each other only for the structure of the approximator chosen
\cite{SP96a}. Most of the papers deal with indirect adaptive control,
trying first to identify the dynamics of the systems and eventually
generating a control input according to the certainty equivalence
principle (with some modification to add robustness to the control
law), whereas very few authors (\eg in \cite{SP96a, RC95}) use the
direct approach, in which the controller {\em directly} generates the
control input to guarantee stability.

Plants whose dynamics can be expressed in the so called ``strict
feedback form'' have been considered, and techniques like backstepping
and adaptive backstepping \cite{K:K:K:1995} have emerged for their
control. The papers \cite{P96, PolMear97} present an extension of the
tuning functions approach in which the nonlinearities of the strict
feedback system are not assumed to be parametric uncertainties, but
rather completely unknown nonlinearities to be approximated on-line
with nonlinearly parameterized function approximators. Both the
adaptive methods in \cite{K:K:K:1995} and in \cite{P96, PolMear97}
attempt to approximate the dynamics of the plant on-line, so they may
be classified as indirect adaptive schemes.

In this paper, we have combined an extension of the class of strict
feedback systems considered in \cite{P96, PolMear97} with the concept
of a dynamic structure that depends on time, so as to propose a class
of nonlinear systems with a time-varying structure, for which we
develop a \emph{direct} adaptive control approach. This class of
systems is a generalization of the class of strict feedback systems
traditionally considered in the literature. Moreover, the direct
adaptive control developed here is, to our knowledge, the first of its
kind in this context, and it presents several advantages with respect
to indirect adaptive methods, including the fact that it needs less
plant information to be implemented.

\section{Direct Adaptive Control}
\label{dac_sec}

Consider the class of continuous time nonlinear systems given by
\begin{align} \label{dac_plant}
  \dot{x}_i &= \sum_{j=1}^R \rho_j(v) \left( \phi_i^j(X_i) +
  \psi_i^j(X_i) x_{i+1} \right) \notag \\
  \dot{x}_n &= \sum_{j=1}^R \rho_j(v) \left( \phi_n^j(X_n) +
  \psi_n^j(X_n) u \right)
\end{align}
where $i=1,2,\dots,n-1$, $X_i = [x_1, \dots, x_i]^\top$, and $X_n \in
\rreal^n$ is the state vector, which we assume measurable, and $u \in
\rreal$ is the control input. The variable $v \in \rreal^q$ may be an
additional input or a possibly exogenous ``scheduling variable.'' We
assume that $v$ and its derivatives up to and including the
$(n-1)^{th}$ one are bounded and available for measurement, which may
imply that $v$ is given by an external dynamical system. The functions
$\rho_j$, $j=1,\dots,R$ may be considered to be ``interpolating
functions'' that produce the time-varying structural nature of system
\eqref{dac_plant}, since they combine $R$ systems in strict feedback
form (given by the $\phi_i^j$ and $\psi_i^j$ functions, $i=1,\dots,n$,
$j=1\dots,R$) and the combination depends on time through the variable
$v$ (thereby, the dynamics of the plant may be different at each time
point depending on the scheduling variable).  Here, we assume that the
functions $\rho_j$ are $n$ times continuously differentiable, and that
they satisfy, for all $v \in \rreal^q$, $ \sum_{j=1}^R \rho_j(v) <
\infty$ and $\left|\frac{\partial^i \rho_j(v)}{\partial v^i} \right| <
\infty$. Denote for convenience $\phi_i^c(X_i, v) = \sum_{j=1}^R
\rho_j(v) \phi_i^j(X_i)$ and $\psi_i^c(X_i, v) = \sum_{j=1}^R
\rho_j(v) \psi_i^j(X_i).$ We will assume that $\phi_i^c$ and
$\psi_i^c$ are sufficiently smooth in their arguments, and that they
satisfy, for all $X_i \in \rreal^i$ and $v \in \rreal^q$,
$i=1,\dots,n$, $\phi_i^c(0, v) = 0$ and $\psi_i^c(X_i, v) \neq 0.$

Here, we will develop a direct adaptive control method for the class of
systems \eqref{dac_plant}. We assume that the interpolation functions $\rho_j$
are known, but the functions $\phi_i^j$ and $\psi_i^j$ (which constitute the
underlying time-varying dynamics of the system) are unknown. In an indirect
adaptive methodology one would attempt to identify the unknown functions and
then construct a stabilizing control law based on the approximations to the
plant dynamics. Here, however, we will postulate the existence of an ideal
control law (based on the assumption that the plant belongs to the class of
systems \eqref{dac_plant}) which possesses some desired stabilizing
properties, and then we devise adaptation laws that attempt to approximate the
ideal control equation.  This approximation will be performed within a compact
set $\mathcal{S}_{x_n} \subset \rreal^n$ of arbitrary size which contains the
origin. In this manner, the results obtained are semi-global, in the sense
that they are valid as long as the state remains within $\mathcal{S}_{x_n}$,
but this set can be made as large as desired by the designer. In particular,
with enough plant information it can be made large enough that the state never
exits it, since, as will be shown a bound can be placed on the state
transient. Furthermore, as will be indicated below, the stability can be made
global by using bounding control terms.

For each vector $X_i$ we will assume the existence of a compact set
$\mathcal{S}_{x_i} \subset \rreal^i$ specified by the designer. We
will consider trajectories within the compact sets $\mathcal{S}_{x_i}$,
$i=1,\dots,n$, where the sets are constructed such that $\mathcal{S}_{x_i}
\subset \mathcal{S}_{x_{i+1}}$, for $i=1, \dots, n-1$. We assume the existence
of bounds $\underline{\psi}_i^c$, $\bar{\psi}_i^c \in \rreal$, and
$\psi_{i_d}^c \in \rreal$, $i=1,\dots,n$ (\emph{not necessarily known}), such
that for all $v \in \rreal^q$ and $X_i \in \mathcal{S}_{x_i}$, $i=1, \dots,
n$,
\begin{align} \label{dac_psi_ass}
  0 &< \underline{\psi}_i^c \leq \psi_i^c(X_i,v) \leq \bar{\psi}_i^c <
  \infty \notag \\
  \left|\dot{\psi}_i^c\right| &= \left| \sum_{j=1}^R \left(
  \frac{\partial \rho_j(v)}{\partial v} \dot{v} \psi_i^j(X_i) +
  \rho_j(v) \frac{\partial \psi_i^j(X_i)}{\partial X_i} \dot{X}_i
  \right) \right| \leq \psi_{i_d}^c.
\end{align}
This assumption implies that the affine terms in the plant dynamics have a
bounded gain and a bounded rate of change. Since the functions $\psi_i^c$ are
assumed continuous, they are therefore bounded within $\mathcal{S}_{x_i}$.
Similarly, note that even though the term $|\dot{X}_i|$ may not necessarily be
globally bounded, it will have a constant bound within $\mathcal{S}_{x_i}$ due
to the continuity assumptions we make. Therefore, assumption
\eqref{dac_psi_ass} will always be satisfied within $\mathcal{S}_{x_n}$.
Moreover, in the simplest of cases, the first part of assumption
\eqref{dac_psi_ass} is satisfied globally when the functions $\psi_i^j$ are
constant or sector bounded for all $X_i \in \rreal^i$.

The class of plants \eqref{dac_plant} is, to our knowledge, the most
general class of systems considered so far within the context of
adaptive control based on backstepping. In particular, in both
\cite{K:K:K:1995} and \cite{P96, PolMear97}, which are indirect
adaptive approaches, the input functions $\psi_i^j$ are assumed to be
constant for $i=1,\dots,n$. This assumption allows the authors of
those works to perform a simpler stability analysis, which becomes
more complex in the general case \cite{OP1999e}. Also, the addition of
the interpolation functions $\rho_j$, $j=1,\dots,R$, extends the class
of strict feedback systems to one including systems with a
time-varying structure \cite{OP2000a}, as well as systems falling in
the domain of gain scheduling (where the plant dynamics are identified
at different operating points and then interpolated between using a
scheduling variable). Note that if we let $R=1$ and $\rho_1(v)=1$ for
all $v$, together with $\psi_i^c=1$, $i=1, \dots, n$, we have the
particular case considered in \cite{P96, PolMear97}.

The direct approach presented here has several advantages with respect to
indirect approaches such as in \cite{K:K:K:1995, P96, PolMear97}.  In
particular, bounds on the input functions $\psi_i^j$ are only assumed to
exist, but need \emph{neither to be known nor to be estimated}. This is
because the ideal law is formulated so that there is not an explicit need to
include information about the bounds in the actual control law.  Moreover,
although assumption \eqref{dac_psi_ass} appears to be more restrictive than
what is needed in the indirect adaptive case, it is in fact not so due to the
fact that the stability results are semi-global (\ie since we are operating
within the compact sets $\mathcal{S}_{x_n}$, continuity of the affine terms
automatically implies the satisfaction of the second part of assumption
\eqref{dac_psi_ass}).

\subsection{Direct Adaptive Control Theorem}
\label{dac_proof}

Next, we state our main result and then show its proof\footnote{We
  will generally omit the arguments of functions for brevity.}. For
convenience, we use the notation $\nu_i=[v, \dot{v}, \dots, v^{(i-1)}]
\in \rreal^{q \times i}$, $i=1,\dots,n$.

\begin{theoremm} \label{dac_theorem}
Consider system \eqref{dac_plant} with the state vector $X_n$ measurable and
the scheduling matrix $\nu_{n-1}$ measurable and bounded, together with
the above stated assumptions on $\phi_i^c$, $\psi_i^c$ and $\rho_j$, and
\eqref{dac_psi_ass}. Assume also that $\nu_i(0) \in \mathcal{S}_{v_i} \subset
\rreal^{q \times i}$, $X_i(0) \in \mathcal{S}_{x_i} \subset \rreal^i$,
$i=1,\dots,n$, where $\mathcal{S}_{v_i}$ and $\mathcal{S}_{x_i}$ are compact
sets specified by the designer, and large enough that $\nu_i$ and $X_i$ do not
exit them. Consider the diffeomorphism $z_1 = x_1$, $z_i = x_i -
\hat{\alpha}_{i-1} - \alpha_{i-1}^s$, $i=2,\dots,n$, with
$\hat{\alpha}_i(X_i, \nu_i) = \sum_{j=1}^R \rho_j(v)
\hat{\theta}_{\alpha_i^j}^\top \zeta_{\alpha_i^j}(X_i, \nu_i)$ and
$\alpha_i^s(z_i, z_{i-1}) = -k_i z_i - z_{i-1}$,
with $k_i > 0$ and $z_0 = 0$. Assume the functions
$\zeta_{\alpha_i^j}(X_i, \nu_i)$ to be at least $n-i$ times
continuously differentiable, and to satisfy, for $i=1,\dots,n$,
$j=1,\dots,R$,
\begin{equation} \label{dac_zeta_ass}
  \left| \frac{\partial^{n-i} \zeta_{\alpha_i^j}}{\partial [X_i,
  \nu_i]^{n-i}} \right| < \infty.
\end{equation}
Consider the adaptation laws for the parameter vectors
$\hat{\theta}_{\alpha_i^j} \in \rreal^{N_{\alpha_i^j}}$,
$N_{\alpha_i^j} \in \nats$, $\dot{\hat{\theta}}_{\alpha_i^j} = -\rho_j
\gamma_{\alpha_i^j} \zeta_{\alpha_i^j} z_i - \sigma_{\alpha_i^j}
\hat{\theta}_{\alpha_i^j}$, where $\gamma_{\alpha_i^j}>0$,
$\sigma_{\alpha_i^j}>0$, $i=1,\dots,n$, $j=1,\dots,R$ are design
parameters. Then, the control law $u = \hat{\alpha}_n + \alpha_n^s$
guarantees boundedness of all signals and convergence of the states to
the residual set
\begin{equation} \label{dac_residual_set}
  \mathcal{D}_d = \left\{ X_n \in \Re^n: \sum_{i=1}^n z_i^2 \leq
  \frac{2 \underline{\psi}_m W_d}{\beta_d} \right\}.
\end{equation}
where $\underline{\psi}_m = \min_{1 \leq i \leq n} \bar{\psi}_i^c$, $\beta_d$
is a constant, and $W_d$ measures approximation errors and ideal parameter
sizes, and its magnitude can be reduced through the choice of the design
constants $k_i$, $\gamma_{\alpha_i^j}$ and $\sigma_{\alpha_i^j}$.
\end{theoremm}

\begin{proof}
The proof requires $n$ steps, and is performed inductively. First, let $z_1 =
x_1$, and $z_2 = x_2 - \hat{\alpha}_1 - \alpha_1^s$, where $\hat{\alpha}_1$ is
the approximation to an ideal signal $\alpha_1^*$ (``ideal'' in the sense that
if we had $\hat{\alpha}_1 = \alpha_1^*$ we would have a globally
asymptotically stable closed loop without need for the stabilizing term
$\alpha_1^s$), and $\alpha_1^s$ will be given below. Let $c_1>0$ be a constant
such that $c_1 > \frac{\psi_{1_d}^c}{2 \underline{\psi}_1^c}$, and
$  \alpha_1^*(x_1,v) = \frac{1}{\psi_1^c} \left( -\phi_1^c - c_1 z_1
  \right).
$
Since the ideal control $\alpha_1^*$ is smooth, it may be approximated with
arbitrary accuracy for $v$ and $x_1$ within the compact sets
$\mathcal{S}_{v_1} \subset \rreal^q$ and $\mathcal{S}_{x_1} \subset \rreal$,
respectively, as long as the size of the approximator can be made arbitrarily
large.

For approximators of finite size let $ \alpha_1^*(x_1,v) = \sum_{j=1}^R
\rho_j(v) \theta_{\alpha_1^j}^{*^\top} \zeta_{\alpha_1^j}(v, x_1) +
\delta_{\alpha_1}(v, x_1), $ where the parameter vectors
$\theta_{\alpha_1^j}^* \in \rreal^{N_{\alpha_1^j}}$, $N_{\alpha_1^j} \in
\nats$, are optimum in the sense that they minimize the representation error
$\delta_{\alpha_1}$ over the set $\mathcal{S}_{x_1} \times \mathcal{S}_{v_1}$
and suitable compact parameter spaces $\Omega_{\alpha_1^j}$, and
$\zeta_{\alpha_1^j}(x_1, v)$ are defined via the choice of the approximator
structure (see \cite{OP2000b} for an example of a choice for
$\zeta_{\alpha_i^j}$). The parameter sets $\Omega_{\alpha_1^j}$ are simply
mathematical artifacts.  As a result of the stability proof the approximator
parameters are bounded using the adaptation laws in Theorem \ref{dac_theorem},
so $\Omega_{\alpha_1^j}$ does not need to be defined explicitly, and no
parameter projection (or any other ``artificial'' means of keeping the
parameters bounded) is required. The representation error $\delta_{\alpha_1}$
arises because the sizes $N_{\alpha_i^j}$ are finite, but it may be made
arbitrarily small within $\mathcal{S}_{x_1} \times \mathcal{S}_{v_1}$ by
increasing $N_{\alpha_i^j}$ (\ie we assume the chosen approximator structures
possess the ``universal approximation property''). In this way, there exists a
constant bound $d_{\alpha_1}>0$ such that $|\delta_{\alpha_1}| \leq
d_{\alpha_1} < \infty$. To make the proof logically consistent, however, we
need to assume that some knowledge about this bound and a bound on
$\theta_{\alpha_1^j}^*$ are available (since in this case it becomes possible
to guarantee a priori that $\mathcal{S}_{x_1} \times \mathcal{S}_{v_1}$ is
large enough). However, in practice some amount of redesign may be required,
since these bounds are typically guessed by the designer

Let $\Phi_{\alpha_1^j} = \hat{\theta}_{\alpha_1^j} - \theta_{\alpha_1^j}^*$
denote the parameter error, and approximate $\alpha_1^*$ with
$\hat{\alpha}_1(x_1, v, \hat{\theta}_{\alpha_1^j}; j=1,\dots,R) = \sum_{j=1}^R
\rho_j(v) \hat{\theta}_{\alpha_1^j}^\top \zeta_{\alpha_1^j}(x_1, v).$ Hence,
we have a linear in the parameters approximator with parameter vectors
$\hat{\theta}_{\alpha_1^j}$. Note that the structural dependence on time of
system \eqref{dac_plant} is reflected in the controller, because
$\hat{\alpha}_1$ can be viewed as using the functions $\rho_j(v)$ to
interpolate between ``local'' controllers of the form
$\hat{\theta}_{\alpha_1^j}^\top \zeta_{\alpha_1^j}(x_1, v)$, respectively.
Notice that since the functions $\rho_j$ are assumed continuous and $v$
bounded, the signal $\hat{\alpha}_1$ is well defined for all $v \in
\mathcal{S}_{v_1}$.

Consider the dynamics of the transformed state, $\dot{z}_1 = \phi_1^c
+ \psi_1^c(z_2 + \hat{\alpha}_1 + \alpha_1^s) + \psi_1^c(\alpha_1^* -
\alpha_1^*) = -c_1 z_1 + \psi_1^c z_2 + \psi_1^c(\hat{\alpha}_1 -
\alpha_1^*) + \psi_1^c \alpha_1^s = -c_1 z_1 + \psi_1^c z_2 + \psi_1^c
\left( \sum_{j=1}^R \rho_j \Phi_{\alpha_1^j}^\top \zeta_{\alpha_1^j} -
  \delta_{\alpha_1^j} \right) + \psi_1^c \alpha_1^s.$ Let $V_1 =
\frac{1}{2 \psi_1^c} z_1^2 + \frac{1}{2} \sum_{j=1}^R
\frac{\Phi_{\alpha_1^j}^\top \Phi_{\alpha_1^j}}{\gamma_{\alpha_1^j}}$,
and examine its derivative, $ \dot{V}_1 = \frac{2 \psi_1^c (2 z_1
  \dot{z}_1) - 2 z_1^2 \dot{\psi}_1^c}{4 \psi_1^{c^2}} + \sum_{j=1}^R
\frac{\Phi_{\alpha_1^j}^\top
  \dot{\Phi}_{\alpha_1^j}}{\gamma_{\alpha_1^j}}.$ Using the expression
for $\dot{z}_1$, $ \dot{V}_1 = -\frac{c_1 z_1^2}{\psi_1^c} + z_1 z_2 +
z_1 \sum_{j=1}^R \rho_j \Phi_{\alpha_1^j}^\top \zeta_{\alpha_1^j} -
z_1 \delta_{\alpha_1^j} + z_1 \alpha_1^s - \frac{1}{2} z_1^2
\frac{\dot{\psi}_1^c}{\psi_1^{c^2}} + \sum_{j=1}^R
\frac{\Phi_{\alpha_1^j}^\top
  \dot{\Phi}_{\alpha_1^j}}{\gamma_{\alpha_1^j}}. $ Choose the
adaptation law $ \dot{\hat{\theta}}_{\alpha_1^j} =
\dot{\Phi}_{\alpha_1^j} = -\rho_j \gamma_{\alpha_1^j}
\zeta_{\alpha_1^j} z_1 - \sigma_{\alpha_1^j}
\hat{\theta}_{\alpha_1^j}, $ with design constants
$\gamma_{\alpha_1^j} > 0$, $\sigma_{\alpha_1^j} > 0$, $j=1,\dots,R$
(we think of $\sigma_{\alpha_1^j} \hat{\theta}_{\alpha_1^j}$ as a
``leakage term''). Also, note that for any constant $k_1>0$, $ -z_1
\delta_{\alpha_1^j} \leq |z_1| d_{\alpha_1} \leq k_1 z_1^2 +
\frac{d_{\alpha_1}^2}{4 k_1}. $ We pick $ \alpha_1^s = -k_1 z_1. $

Notice also that, completing squares, $ -\Phi_{\alpha_1^j}^\top
\hat{\theta}_{\alpha_1^j} = -\Phi_{\alpha_1^j}^\top (\Phi_{\alpha_1^j}
+ \theta_{\alpha_1^j}^*) \leq - \frac{| \Phi_{\alpha_1^j}|^2}{2} +
\frac{| \theta_{\alpha_1^j}^* |^2}{2}. $ Finally, observe that $
-\frac{z_1^2}{\psi_1^c} \left(c_1 + \frac{\dot{\psi}_1^c}{2 \psi_1^c}
\right) \leq -\frac{z_1^2}{\psi_1^c} \left(c_1 - \frac{\psi_{1_d}^c}{2
    \underline{\psi}_1^c} \right) \leq -\frac{\bar{c}_1
  z_1^2}{\bar{\psi}_1^c}, $ with $\bar{c}_1 = c_1 -
\frac{\psi_{1_d}^c}{2\underline{\psi}_1^c} > 0$. Then, we obtain $
\dot{V}_1 \leq -\frac{\bar{c}_1 z_1^2}{\bar{\psi}_1^c} - \frac{1}{2}
\sum_{j=1}^R \sigma_{\alpha_1^j} \frac{| \Phi_{\alpha_1^j}
  |^2}{\gamma_{\alpha_1^j}} + z_1 z_2 + \frac{d_{\alpha_1}^2}{4 k_1} +
\frac{1}{2} \sum_{j=1}^R \sigma_{\alpha_1^j}
\frac{\theta_{\alpha_1^j}^*}{\gamma_{\alpha_1^j}}. $ This completes
the first step of the proof.

We may continue in this manner up to the $n^{th}$ step\footnote{We
  omit intermediate steps for brevity.}, where we have $z_n = x_n -
\hat{\alpha}_{n-1} - \alpha_{n-1}^s$, with $\hat{\alpha}_{n-1}$ and
$\alpha_{n-1}^s$ defined as in Theorem \ref{dac_theorem}. Consider the
ideal signal $ \alpha_n^*(X_n, \nu_n) = \frac{1}{\psi_n^c} \left(
  \phi_n^c - c_n z_n + \dot{\hat{\alpha}}_{n-1} + \dot{\alpha}_{n-1}^s
\right) $ with $c_n > \frac{\psi_{n_d}^c}{2 \underline{\psi}_n^c}$.
Notice that, even though the terms
$\dot{\hat{\theta}}_{\alpha_{n-1}^j}$ appear in $\alpha_n^*$ through
the partial derivatives in $\dot{\hat{\alpha}}_{n-1}$,
$\hat{\theta}_{\alpha_{n-1}^j}$ does not need to be an input to
$\alpha_n^*$, since the resulting product of the partial derivatives
and $\dot{\hat{\theta}}_{\alpha_{n-1}^j}$ can be expressed in terms of
$z_1,\dots,z_{n-1}$, $v$ and $\sigma_{\alpha_{n-1}^j}
\hat{\alpha}_{n-1}$. To simplify the notation, however, we will omit
the dependencies on inputs other than $X_i$ and $\nu_i$, but bearing
in mind that, when implementing this method, more inputs may be
required to satisfy the proof. Also, note that by assumption
\eqref{dac_zeta_ass}, $|\alpha_n^*| < \infty$ for bounded arguments.
Therefore, we may represent $\alpha_n^*$ with $ \alpha_n^*(X_n, \nu_n)
= \sum_{j=1}^R \rho_j(v) \theta_{\alpha_n^j}^{*^\top}
\zeta_{\alpha_n^j}(X_n, \nu_n) + \delta_{\alpha_n}(X_n, \nu_n) $ for
$X_n \in \mathcal{S}_{x_n} \subset \rreal^n$ and $\nu_n \in
\mathcal{S}_{v_n} \subset \rreal^{q \times n}$. The parameter vector
$\theta_{\alpha_n^j}^* \in \rreal^{N_{\alpha_n^j}}$, $N_{\alpha_n^j}
\in \nats$ is an optimum within a compact parameter set
$\Omega_{\alpha_n}$, in a sense similar to $\theta_{\alpha_1^j}^*$, so
that for $(X_n, \nu_n) \in \mathcal{S}_{x_n} \times
\mathcal{S}_{v_n}$, $|\delta_{\alpha_n}| \leq d_{\alpha_n} < \infty$
for some bound $d_{\alpha_n} > 0$. Let $\Phi_{\alpha_n^j} =
\hat{\theta}_{\alpha_n^j} - \theta_{\alpha_n^j}^*$, and consider the
approximation $\hat{\alpha}_n$ as given in Theorem \ref{dac_theorem}.
The control law $u = \hat{\alpha}_n + \alpha_n^s$ yields $ \dot{z}_n =
\phi_n^c + \psi_n^c(\hat{\alpha}_n + \alpha_n^s) -
\dot{\hat{\alpha}}_{n-1} - \dot{\alpha}_{n-1}^s + \psi_n^c(\alpha_n^*
- \alpha_n^*) = -c_n z_n + \psi_n^c \left( \sum_{j=1}^R \rho_j(v)
  \Phi_{\alpha_n^j}^\top \zeta_{\alpha_n^j} - \delta_{\alpha_n}
\right) + \psi_n^c \alpha_n^s. $ Choose the Lyapunov function
candidate $ V = V_{n-1} + \frac{1}{2 \psi_n^c} z_n^2 + \frac{1}{2}
\sum_{j=1}^R \frac{\Phi_{\alpha_n^j}^\top
  \Phi_{\alpha_n^j}}{\gamma_{\alpha_n^j}} $ and examine its
derivative, $ \dot{V} = \dot{V}_{n-1} - \frac{c_n z_n^2}{\psi_n^c} +
z_n \sum_{j=1}^R \rho_j(v)\Phi_{\alpha_n^j}^\top \zeta_{\alpha_n^j} -
z_n \delta_{\alpha_n} + z_n \alpha_n^s - \frac{1}{2} z_n^2
\frac{\dot{\psi}_n^c}{\psi_n^{c^2}} +
\sum_{j=1}^R\frac{\Phi_{\alpha_n^j}^\top
  \dot{\Phi}_{\alpha_n^j}}{\gamma_{\alpha_n^j}}$. One can show
inductively that $ \dot{V}_{n-1} \leq -\sum_{i=1}^{n-1}
\frac{\bar{c}_i z_i^2}{\bar{\psi}_i^c} - \frac{1}{2} \sum_{i=1}^{n-1}
\sum_{j=1}^R \sigma_{\alpha_i^j}
\frac{|\Phi_{\alpha_i^j}|^2}{\gamma_{\alpha_i^j}} + z_{n-1} z_n +
\sum_{i=1}^{n-1} \frac{d_{\alpha_i}^2}{4 k_i} + \frac{1}{2}
\sum_{i=1}^{n-1} \sum_{j=1}^R \sigma_{\alpha_i^j}
\frac{|\theta_{\alpha_i^j}^*|^2}{\gamma_{\alpha_i^j}} $ with constants
$\bar{c}_i = c_i - \frac{\psi_{i_d}}{2 \underline{\psi}_i^c} > 0$,
$i=1,\dots,n$. The choice of adaptation laws for $\theta_{\alpha_n^j}$
and of $\alpha_n^s$ in Theorem \ref{dac_theorem}, together with the
observations that $ -\frac{\sigma_{\alpha_n^j}}{\gamma_{\alpha_n^j}}
\Phi_{\alpha_n^j}^\top \hat{\theta}_{\alpha_n^j} \leq
-\frac{\sigma_{\alpha_n^j}}{\gamma_{\alpha_n^j}}
\frac{|\Phi_{\alpha_n^j}|^2}{2} +
\frac{\sigma_{\alpha_n^j}}{\gamma_{\alpha_n^j}}
\frac{|\theta_{\alpha_n^j}^*|^2}{2}$, $ - z_n \delta_{\alpha_n^j} \leq
k_n z_n^2 + \frac{d_{\alpha_n}}{4 k_n}$, with $k_n > 0$ and $
-\frac{z_n^2}{\psi_n^c} \left(c_n + \frac{\dot{\psi}_n^c}{2 \psi_n^c}
\right) \leq -\frac{\bar{c}_n z_n^2}{\bar{\psi}_n^c} $ imply
\begin{equation}
  \dot{V} \leq -\sum_{i=1}^n \frac{\bar{c}_i z_i^2}{\bar{\psi}_i^c} -
  \frac{1}{2} \sum_{i=1}^n \sum_{j=1}^R \sigma_{\alpha_i^j}
  \frac{|\Phi_{\alpha_i^j}|^2}{\gamma_{\alpha_i^j}} + W_d,
\end{equation}
where $W_d$ contains the combined effects of representation errors and ideal
parameter sizes, and is given by
$
  W_d = \sum_{i=1}^n \frac{d_{\alpha_i}^2}{4 k_i} + \frac{1}{2}
  \sum_{i=1}^n \sum_{j=1}^R \sigma_{\alpha_i^j}
  \frac{|\theta_{\alpha_i^j}^*|^2}{\gamma_{\alpha_i^j}}.
$
Note that if
$
  \sum_{i=1}^n \frac{\bar{c}_i z_i^2}{\bar{\psi}_i^c} \geq W_d$ or
$  \frac{1}{2} \sum_{i=1}^n \sum_{j=1}^R \sigma_{\alpha_i^j}
  \frac{|\Phi_{\alpha_i^j}|^2}{\gamma_{\alpha_i^j}} \geq W_d
$,
then we have $\dot{V} \leq 0$. Furthermore, letting $\underline{\psi}_m =
\min_{1 \leq i \leq n} (\underline{\psi}_i^c)$, $\bar{\psi}_m = \max_{1 \leq i
  \leq n} (\bar{\psi}_i^c)$, and defining
$
  \bar{c}_0 = \min_{1\leq i\leq n} (\bar{c}_i)$,
$  \psi_m = \frac{\underline{\psi}_m}{\bar{\psi}_m}$ and
$  \sigma_0 = \min_{1 \leq i \leq n, 1 \leq j \leq R}
  \left( \sigma_{\alpha_i^j} \right)$
we have
$
  -\sum_{i=1}^n \frac{\bar{c}_i z_i^2}{\bar{\psi}_i^c} \leq
  -\bar{c}_0 \sum_{i=1}^n \frac{z_i^2}{\bar{\psi}_i^c} = -\bar{c}_0
  \sum_{i=1}^n \frac{z_i^2}{\psi_i^c} \frac{\psi_i^c}{\bar{\psi}_i^c}
  \leq -\bar{c}_0 \sum_{i=1}^n \frac{z_i^2}{\psi_i^c}
  \frac{\underline{\psi}_i^c}{\bar{\psi}_i^c} \leq -\bar{c}_0
  \psi_m \sum_{i=1}^n \frac{z_i^2}{\psi_i^c}$ and
$  -\frac{1}{2} \sum_{i=1}^n \sum_{j=1}^R \sigma_{\alpha_i^j}
  \frac{|\Phi_{\alpha_i^j}|^2}{\gamma_{\alpha_i^j}} \leq -\sigma_0
  \frac{1}{2} \sum_{i=1}^n \sum_{j=1}^R
  \frac{|\Phi_{\alpha_i^j}|^2}{\gamma_{\alpha_i^j}}.
$
Then, letting $\beta_d = \min (2 \bar{c}_0 \psi_m,
\sigma_0)$, we have that if
\begin{equation}
  V = \frac{1}{2} \sum_{i=1}^n \frac{z_i^2}{\psi_i^c} + \frac{1}{2}
  \sum_{i=1}^n \sum_{j=1}^R
  \frac{|\Phi_{\alpha_i^j}|^2}{\gamma_{\alpha_i^j}} \geq V_0
\end{equation}
with $V_0 = \frac{W_d}{\beta_d}$, then $\dot{V} \leq 0$ and all
signals in the closed loop are bounded. Furthermore, we have $\dot{V}
\leq -\beta_d V + W_d$, which implies that $ 0 \leq V(t) \leq
\frac{W_d}{\beta_d} + \left( V(0) - \frac{W_d}{\beta_d} \right)
e^{-\beta_d t} $ so that both the transformed states and the parameter
error vectors converge to a bounded set. Finally, we conclude from the
upper bound on $V(t)$ that the state vector $X_n$ converges to the
residual set \eqref{dac_residual_set}.
\end{proof}

\begin{remarkk} \label{dac_rem_bounding_control}
  The representation error bounds and the size of the ideal parameter vectors
  are assumed known, since they affect the size of the residual set to which
  the states converge. It is possible to augment the direct adaptive algorithm
  with ``auto-tuning'' capabilities (similar to \cite{PolMear97}), which would
  relax the need for these bounds.
  
  Furthermore, note that the stability result of Theorem \ref{dac_theorem} is
  semi-global, in the sense that it is valid within the compact sets
  $\mathcal{S}_{v_i}$ and $\mathcal{S}_{x_i}$, $i=1,\dots,n$, which can be
  made arbitrarily large. The stability result may be made global by adding a
  high gain bounding control term to the control law. Such a term may be
  particularly useful when, due to a complete lack of a priori knowledge, the
  control designer is unable to guarantee that the compact sets
  $\mathcal{S}_{x_i}$, $i=1,\dots,n$, are large enough so that the state will
  not exit them before the controller has time to bring the state inside
  $\mathcal{D}_d$; moreover, it may also happen that due to a poor design and
  poor system knowledge, $\mathcal{D}_d$ is not contained in
  $\mathcal{S}_{x_n}$. In this case, too, bounding control terms may be
  helpful until the design is refined and improved. However, using bounding
  control requires explicit knowledge of functional upper bounds of
  $|\psi_i^c(v, X_i)|$, and also of the lower bounds $\underline{\psi}_i^c$,
  $i=1,\dots,n$, whose knowledge we do not mandate in Theorem
  \ref{dac_theorem}. Bounding terms may be added to the diffeomorphism in
  Theorem \ref{dac_theorem}, but we do not present the analysis since it is
  similar to the one we present here and it is algebraically tedious; we
  simply note, though, that the bounding terms have to be smooth (because they
  need to be differentiable), so they need to be defined in terms of smooth
  approximations to the sign, saturation and absolute value functions that are
  typically used in this approach.
\end{remarkk}

\begin{remarkk} \label{dac_setc_i}
  If the bounds $\underline{\psi}_i^c$, $\bar{\psi}_i^c$ and $\psi_{i_d}^c$
  are known, it becomes possible for the designer to directly set the
  constants $c_i$ in the control law. Notice that with knowledge of these
  bounds, the term $\underline{\psi}_m$ is also known, and we can pick
  constants $c_i$ such that $c_i > \frac{\psi_{i_d}^c}{2
    \underline{\psi}_i^c}$. Define the auxiliary functions $\eta_i = c_i z_i$.
  We may explicitly set the constant $c_i$ in $\alpha_i^*$ if we let $\eta_i$
  be an input to the $i^{th}$ approximator structure, \ie if we let
  $\alpha_i^*(X_i, \nu_i, \dot{X}_{r_i}, \eta_i) = \sum_{j=1}^R \rho_j(v)
  \theta_{\alpha_i^j}^{*^\top} \zeta_{\alpha_i^j}(X_i, \nu_i, \dot{X}_{r_i},
  \eta_i) + \delta_{\alpha_i}$.  Then, the approximators used in the control
  procedure are given by $\hat{\alpha}_i(X_i, \nu_i, \dot{X}_{r_i}, \eta_i) =
  \sum_{j=1}^R \rho_j(v) \hat{\theta}_{\alpha_i^j}^\top
  \zeta_{\alpha_i^j}(X_i, \nu_i, \dot{X}_{r_i}, \eta_i)$ and the stability
  analysis can be carried out as expected.
\end{remarkk}

\subsection{Performance Analysis: $\Lnorm_\infty$ Bounds and Transient
Design}
\label{dac_transient}

The stability result of Theorem \ref{dac_theorem} is useful in that it
indicates conditions to obtain a stable closed-loop behavior for a
plant belonging to the class given by \eqref{dac_plant}. However, it
is not immediately clear how to choose the several design constants to
improve the control performance. Here we concentrate on the tracking
problem, and present design guidelines with respect to an
$\Lnorm_\infty$ bound on the tracking error. We are interested in
having $x_1$ track the reference model state $x_{r_1}$ of the
reference model $\dot{x}_{r_i} = x_{r_{i+1}}, \text { } i=1, 2, \dots,
n-1$, $\dot{x}_{r_n} = f_r(X_{r_n}, r)$ with bounded reference input
$r(t) \in \rreal$. Now, we need to use the diffeomorphism $z_1 = x_1 -
x_{r_1}$, $z_i = x_i - \hat{\alpha}_{i-1} - \alpha_{i-1}^s$,
$i=2,\dots,n$ with $\alpha_1^*(x_1, v, \dot{x}_{r_1}) =
\frac{1}{\psi_1^c} \left( -\phi_1^c - c_1 z_1 + x_{r_2} \right)$ and
$\alpha_i^*(X_i, \nu_i, \dot{X}_{r_i}) = \frac{1}{\psi_i^c} \left(
  -\phi_i^c - c_i z_i + \dot{\hat{\alpha}}_i + \dot{\alpha}_i^s
\right)$ for $i=2, \dots, n$. The stability proof needs to be modified
accordingly, and it can be shown that the tracking error $|x_1 -
x_{r_1}|$ converges to a neighborhood of size $\sqrt{\frac{2
    \underline{\psi}_m W_d}{\beta_d}}$.

From the upper bound on $V(t)$ we can write $ V(t) \leq
\frac{W_d}{\beta_d} + V(0) e^{-\beta_d t}$. From here, it follows that
$ \frac{1}{2} \sum_{i=1}^n \frac{z_i^2(t)}{\psi_i^c(t)} \leq
\frac{W_d}{\beta_d} + \left( \frac{1}{2} \sum_{i=1}^n
  \frac{z_i^2(0)}{\psi_i^c(0)} + \frac{1}{2} \sum_{i=1}^n \sum_{j=1}^R
  \frac{|\Phi_{\alpha_i^j}(0)|^2}{\gamma_{\alpha_i^j}} \right)
e^{-\beta_d t}. $ The terms $z_i(0)$ depend on the design constants in
a complex manner.  For this reason, rather than trying to take them
into account in the design procedure, we follow the trajectory
initialization approach taken in \cite{K:K:K:1995}, which allows the
designer to set $z_i(0)=0$, $i=1, \dots, n$ by an appropriate choice
of the reference model's initial conditions.  In our case, in addition
to the assumption that it is possible to set the initial conditions of
the reference model, we will have to assume certain invertibility
conditions on the approximators.  In particular, since $z_1(0) =
x_1(0) - x_{r_1}(0)$, for $z_1(0) = 0$ we need to set $x_{r_1}(0) =
x_1(0)$.

For the $i^{th}$ transformed state $z_{i}$, $i=2,\dots,n$, $z_i(0) =
x_i(0) - \hat{\alpha}_{i-1}(0) - \alpha_{i-1}^s(0)$. Notice that
$\alpha_{i-1}^s(0) = \alpha_{i-1}^s(z_{i-1}(0), z_{i-2}(0))$, so that
if $z_{i-1}(0)=0$ and $z_{i-2}(0)=0$ we have $\alpha_{i-1}^s(0)=0$. In
particular, notice that this holds for $i=2$. In this case, to set
$z_2(0)=0$ we need to have $\hat{\alpha}_1(x_1(0), v(0), x_{r_2}(0)) =
x_2(0)$. This equation can be solved analytically (or numerically) for
$x_{r_2}(0)$ provided $\left. \frac{\p \hat{\alpha}_1}{\p x_{r_2}}
\right|_{t=0} \neq 0$. This is not an unreasonable condition, since it
depends on the choice of approximator structure the designer makes.
The structure can be chosen so that it satisfies this condition.
Granted this is the case, it clearly holds that $\alpha_2^s(0)=0$, and
the same procedure can be inductively carried out for $i=3,\dots,n$,
with the choices $\hat{\alpha}_{i-1}(X_{i-1}(0), \nu_{i-1}(0),
x_{r_i}(0)) = x_i(0)$.

This procedure yields the simpler bound $\sum_{i=1}^n z_i^2(t) \leq
\frac{2 \underline{\psi}_m W_d}{\beta_d} + \underline{\psi}_m \left(
  \sum_{i=1}^n \sum_{j=1}^R \frac{|\Phi_{\alpha_i^j}(0)|^2}
  {\gamma_{\alpha_i^j}} \right) e^{-\beta_d t}$. We would like to make
this bound small, so that the transient excursion of the tracking
error is small. Notice that we do not have direct control on the size
of $\beta_d$, since this term depends on the unknown constants $c_i$,
which appear in the ideal signals $\alpha_i^*$. Even though it is not
necessary to be able to set $\beta_d$ to reduce the size of the bound,
it is possible to do so if the bounds $\underline{\psi}_i^c$,
$\bar{\psi}_i^c$ and $\psi_{i_d}^c$ are known.

At this point, it becomes more clear how to choose the constants to
achieve a smaller bound. Recalling the expression of $W_d$, note that,
first, one may want to have $\beta_d > 1$, so that $W_d$ is not made
larger when divided by $\beta_d$, and so that the convergence is
faster. This may be achieved by setting $c_i$ such that $2 \bar{c}_i
\psi_m > 1$ (if enough knowledge is available to do so) and
$\sigma_{\alpha_i^j} > 1$. However, having large $\sigma_{\alpha_i^j}$
makes $W_d$ larger; this can be offset, however, by also choosing the
ratio $\sigma_{\alpha_i^j} / \gamma_{\alpha_i^j} < 1$ or smaller.
Finally, it is clear that making $k_i$ larger reduces the effects of
the representation errors, and therefore makes $W_d$ smaller. Observe
that there is enough design freedom to make $W_{d}$ small and
$\beta_{d}$ large independently of each other.

These simple guidelines may become very useful when performing a real
control design. Moreover, notice that the bound on $\sum_{i=1}^n
z_i^2(t)$ makes it possible to specify the compact sets of the
approximators so that, even throughout the transient, it can be
guaranteed that the states will remain within the compact sets without
the need for a global bounding control term. This has been a recurrent
shortcoming of many on-line function approximation based methods, and
the explicit bound on the transient makes it possible to overcome it.

\section{Conclusions}
\label{adapt_tv_bs_concl}

In this paper we have developed a direct adaptive control method for a class
of uncertain nonlinear systems with a time-varying structure using a Lyapunov
approach to construct the stability proofs. The systems we consider are
composed of a finite number of ``pieces,'' or dynamic subsystems, which are
interpolated by functions that depend on a possibly exogenous scheduling
variable. We assume that each piece is in strict feedback form, and show that
the methods yield stability of all signals in the closed-loop, as well as
convergence of the state vector to a residual set around the equilibrium,
whose size can be set by the choice of several design parameters

We argue that the direct adaptive method presents several advantages over
indirect methods in general, including the need for a smaller amount of
information about the plant and a simpler design. Finally, we provide design
guidelines based on $\Lnorm_\infty$ bounds on the transient and argue that
this bound makes it possible to precisely determine how large the compact sets
for the function approximators should be so that the states do not exit them.

\bibliographystyle{ieeetr}
\bibliography{researchbib,othersbib}

\begin{thebibliography}{10}

\bibitem{Polycarpou91}
M.~M. Polycarpou and P.~A. Ioannou, ``Identification and control of nonlinear
  systems using neural network models: Design and stability analysis,''
  Electrical Engineering -- Systems Report 91-09-01, University of Southern
  California, Sept. 1991.

\bibitem{P96}
M.~M. Polycarpou, ``Stable adaptive neural control scheme for nonlinear
  systems,'' {\em IEEE Transactions on Automatic Control}, vol.~41,
  pp.~447--451, Mar. 1996.

\bibitem{Su94}
C.-Y. Su and Y.~Stepanenko, ``Adaptive control of a class of nonlinear systems
  with fuzzy logic,'' {\em {\IEEETFS}}, vol.~2, pp.~285--294, Nov. 1994.

\bibitem{SP96a}
J.~T. Spooner and K.~M. Passino, ``Stable adaptive control using fuzzy systems
  and neural networks,'' {\em IEEE Transactions in Fuzzy Systems}, vol.~4,
  pp.~339--359, Aug. 1996.

\bibitem{RC95}
G.~A. Rovithakis and M.~A. Christodoulou, ``Direct adaptive regulation of
  unknown nonlinear dynamical systems via dynamic neural networks,'' {\em IEEE
  Transactions on Systems, Man, and Cybernetics}, vol.~25, pp.~1578--1995, Dec.
  1995.

\bibitem{K:K:K:1995}
M.~Krsti\'c, I.~Kanellakopoulos, and P.~Kokotovi\'c, {\em Nonlinear and
  Adaptive Control Design}.
\newblock New York, NY: John Wiley and Sons, 1995.

\bibitem{PolMear97}
M.~M. Polycarpou and M.~J. Mears, ``Stable adaptive tracking of uncertain
  systems using nonlinearly parametrized on-line approximators,'' {\em
  International Journal of Control}, vol.~70, pp.~363--384, May 1998.

\bibitem{OP1999e}
R.~{Ord\'o\~nez} and K.~M. Passino, ``Indirect adaptive control for a class of
  time-varying nonlinear systems,'' {\em Accepted for publication in the
  International Journal of Control}, 2000.

\bibitem{OP2000a}
R.~{{Ord\'o\~nez}} and K.~M. Passino, ``Control of continuous time nonlinear
  systems with a time-varying structure,'' in {\em Proc. of the American
  Control Conf.}, (Chicago, IL), pp.~164--168, June 2000.

\bibitem{OP2000b}
R.~{Ord\'o\~nez} and K.~M. Passino, ``Wing rock regulation with a time-varying
  angle of attack,'' in {\em Proceedings of the Int. Symp. Intelligent
  Control}, (Patras, Greece), pp.~145--150, July 17--19 2000.

\end{thebibliography}

\end{document}